\documentclass[a4paper,12pt]{amsart}
\usepackage{epsf,amscd,amsmath,amssymb,amsxtra,array,sect,color}
\usepackage[xdvi]{graphicx}
\newcommand{\sfit}{\sffamily\slshape}
\newcommand{\R}{\mathbb R}
\newcommand{\Z}{\mathbb Z}

\newcommand{\rnk}{\operatorname{rk}}
\newcommand{\cal}{\mathcal}
\newcommand{\sminus}{\smallsetminus}

\def\fig#1{\vcenter{\hbox{\includegraphics{#1}}}}
\def\sfig#1#2{\vcenter{\hbox{\includegraphics[scale=#2]{#1}}}}

\theoremstyle{plain}
\newtheorem{thm}{\itshape Theorem}

\newtheorem{cor}{\itshape Corollary}

\title[Orientations of Chord Diagrams and Khovanov Homology]{Virtual Links, Orientations of Chord Diagrams 
and Khovanov Homology}

\author{Oleg Viro}
\dedicatory{Uppsala University, Uppsala, Sweden\break
POMI, St.\ Petersburg, Russia}
\address{Department of Mathematics, Uppsala University,
Box 480, S-751 06 Uppsala, Sweden
}
\email{oleg@math.uu.se}

\begin{document}  
\begin{abstract}
By adding or removing appropriate structures to Gauss diagram, one 
can create useful objects related to virtual links. 
In this paper few objects of this kind are studied: 
twisted virtual links generalizing virtual links; 
signed chord diagrams staying halfway between twisted 
virtual links and Kauffman bracket / Khovanov
homology; alternatable virtual links intermediate 
between virtual and classical links. The most profound 
role here belongs to a structure that we dare to call orientation
of chord diagram. Khovanov homology is generalized to oriented 
signed chord diagrams and links in oriented thickened surface such that
the link projection realizes the first Stiefel-Whitney
class of the surface.   
\end{abstract} 
\maketitle

\section{Introduction}\label{s1}

In the middle of nineties Lou Kauffman defined a natural extension
of classical knot theory by replacing classical links with {\it virtual
links.\/} Although introduced in a diagrammatic way, this theory was shown
at its very beginning to be equivalent to study of links in thickened surfaces up to homeomorphism and stabilization by the addition and subtraction of handles 
that do not intersect the link.

Greg Kuperberg \cite{Kup} showed that a representative of a 
virtual link in thickened surface of the minimal genus is defined up to an ambient homeomorphism. Thus, the stabilization can be eliminated at the cost of requiring
minimality of the genus. This provides an opportunity to consider virtual links 100\% traditional objects of 3-dimensional topology: pairs consisting of orientated compact 3-manifold and its oriented closed 1-submanifold up to orientation preserving homeomorphism,  diagrammatic and stabilization free.
 
Therefore diagrammatic technique extended from classical links to virtual ones is extended to links in thickened surfaces. This became
especially interesting due to recent development of link homology
theories built on link diagrams. 

This paper was written in an attempt to analyze the difficulties emerged
when the construction of Khovanov homology is extended to virtual links.
I analyzed few geometric objects closely related to a virtual link and
found that a key role is played by an additional structure, orientation, 
on one of them, a signed chord diagram. Existence of this structure is 
solely responsible for possibility to construct Khovanov complex 
literally as in the classical case. 

Virtual links with orientable chord diagrams appeared in
various occasions: virtual links that admit checkerboard colorable
diagrams, virtual links which can be made alternating by crossing
changes, virtual links with zero-homologous modulo 2 irreducible model,
virtual links with orientable atom, etc. I will call virtual links of
this kind {\sfit alternatable.\/}
 
Construction of Khovanov homology with coefficients in $\Z/2\Z$ 
for any virtual link is quite  straightforward. 
It was presented by V.O.Manturov in \cite{M1}, \cite{M2} and \cite{M3}.
In the same papers Manturov presented also two constructions of
Khovanov homology with integer coefficients, but they do not
generalize the original Khovanov homology of classical links and 
rely on preliminary geometric constructions. 

Manturov's constructions of Khovanov complexes with integer
coefficients  split into geometric
constructions cooking from an arbitrary virtual link alternatable ones, 
and the straightforward construction of Khovanov complex
with integer coefficients which works only for alternatable links.
The first construction is defined
for a framed virtual link, and proceeds by doubling it. The second
one relies on a two-fold covering obtained as restriction to the link of
an orientation covering of a surface, on which the link is naturally 
placed as a zero-homologous modulo 2 curve. This surface with the link
diagram and checkerboard coloring is called the {\sfit atom\/} of the
link. For details see Manturov \cite{M1}, or \cite{M2}, or \cite{M3}.  

The construction of Khovanov complex for alternatable virtual link
does not require all information contained in the virtual link diagram.
One can pass (without any loss) to Gauss diagram and then forget
orientations of its chords. A signed chord diagram obtained in this
way contains everything needed for building of the Khovanov complex.
Alternability of a virtual link is orientability of its chord diagram.

Virtual links admit a generalization to {\sfit twisted virtual
links\/}, which emerge in relation to links in oriented thickennings 
of non-orientable surfaces. Many link invariants, in particular 
Kauffman bracket and Jones polynomial, are extended to twisted virtual links
and links in oriented thickennings of non-orientable surfaces. 

A link in an oriented thickenning of a non-orientable surface which
gives 
rise to oriented signed chord diagram realizes homology class dual to
the first Stiefel-Whitney class of the surface. This is the widest class
of links in oriented thickennings of surfaces for which the classical
construction of Khovanov complex works without any modification over the
integers. 

This created a peculiar situation: the construction of Khovanov complex
with integer coefficients works, say, for non-zero homologous links in 
the real projective space, but does not work for zero-homologous links
in the same space. It is difficult to believe that this is not due just 
to a lack of technique, especially since in the theory of Heegaard-Floer
homology developed by Ozsv\'ath and Szab\'o, which so far was quite similar
to the Khovanov homology, these two classes of links in the projective
space both have homology with integer coefficients. Khovanov,
Rasmussen and Manturov conjectured that there should be a twisted 
version of Khovanov complex which works for all virtual links. 

The main part of the paper starts with yet another introduction to 
virtual knot theory incorporating twisted virtual links. 
Then we review Kauffman bracket and Jones polynomial constructions, 
both for classical and virtual links, and show that they are defined for
any signed chord diagram.
The original part of the paper is devoted to the notion of orientation 
of chord diagram and constructing of Khovanov homology for orientable
signed chord diagrams. 

I am grateful to V.O.Manturov and A.N.Shumakovich for valuable information 
and interesting discussions.  

\subsection{Post-Publication  Remarks} 
After this paper was published, as \cite{V},  V.~O.~Manturov succeeded \cite{M4} in extending Khovanov homology to arbitrary virtual links. This required a deep revision
of the original construction. In particular, in the case of classical link the complex 
defined by Manturov \cite{M4} is not isomorphic to the Khovanov complex (although the complexes are homotopy equivalent).  The new construction also can be modified to
become extendible to arbitrary signed chord diagrams. I am going to devote a separate paper to this subject.

Another generalization of the Khovanov homology is provided by a combination of the categorification of the Kauffman bracket skein module by M.~M.~Asaeda, J.~H.~Przytycki and A.~S.~Sikora  \cite{APS} with the Kuperberg theorem \cite{Kup}. Conjecturally, there should exists a spectral sequence starting at the homology defined by  Asaeda, Przytycki and Sikora \cite{APS} and converging to the homology  defined by Manturov \cite{M4}.
However, the homology defined by Asaeda, Przytycki and Sikora seems to be more appropriate as invariant of a link in a thickened surface: for virtual links it is less explicit due to use of Kuperberg's theorem.

\section{Three Faces of Virtual Knot Theory}\label{s2}

\subsection{Link Diagrams and Gauss Diagrams}\label{s2.1}
To describe graphically a classical link (i.e., a closed 
smooth 1-dimensional submanifold of $\R^3$), 
one takes its generic projection to a plane and 
decorates the image at double points to show over- and under-passes. 
This gives rise to a {\it link diagram:\/}
$$\vcenter{\hbox{\includegraphics[scale=.7]{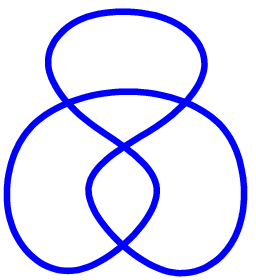}}}  
\qquad \mapsto \qquad 
 \vcenter{\hbox{\includegraphics[scale=.7]{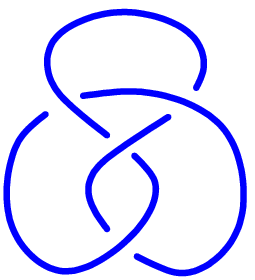}}}$$

Genericity of projection means here that this is an immersion 
with multiple points of multiplicity 2 at most and 
transversality at double points.

A link diagram is a two-dimensional picture of link. 
In many cases one-dimensional picture serves better. 
In particular, it is easier to convert to a combinatorial 
description, used as input data in computer programs.

A one-dimensional picture comes from a parametrization of 
the link:
$$\input{figs/d-fig8-pzed.pstex_t}$$
The source of the parametrization is decorated. 
First, it is oriented.
Second, each over-pass is connected to the corresponding 
under-pass with an arrow:
$$\input{figs/gd-nosigns-fig8.pstex_t}$$
Third, each arrow is equipped with a sign:
$$\input{figs/gd-fig8.pstex_t}$$
The sign is the {\sfit local writhe\/} of the crossing. 
It is $+$ at
a crossing which looks like this:
$\vcenter{\hbox{\includegraphics[scale=.8]{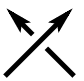}}}$, 
and $-$ at a crossing which looks like that:
$\vcenter{\hbox{\includegraphics[scale=.8]{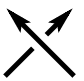}}}$. 
The signs depend on the orientation. 

The result is called a {\sfit Gauss diagram\/} of the link. 
Gauss diagrams were introduced in our joint paper with 
Polyak \cite{PV1}. The corresponding combinatorial objects 
called {\sfit Gauss codes\/} can be traced back to Gauss' 
notebooks. Transition from Gauss codes to their geometric 
counter-parts, Gauss diagrams, is encouraged by geometric 
structures and operations, such as orientations and surgery, 
traditional for geometric objects, but difficult to recognize 
in a combinatorial context.   

A Gauss diagram consists of an oriented 1-manifold (not necessarily 
connected) and chords connecting disjoint pairs of points on the 1-manifold.
Each of the chords is oriented and equipped with a sign. The
1-manifold is called the {\sfit base\/} of the Gauss diagram.

\subsection{From Gauss Diagram to Virtual Link}\label{s2.2}
Not any Gauss diagram can be obtained from a link diagram, but for any Gauss
diagram one can try. 

Take, for example, Gauss diagram $\vcenter{\hbox{\input{figs/GD2VD1col-05.pstex_t}}}$, 
and try to reconstruct the knot. Let us start with crossings, as they
are clearly described up to plane isotopy. Then connect them step by step 
according to the Gauss diagram:
$$ \vcenter{\hbox{\input{figs/GD2VD2col-07.pstex_t}}} \to \ 
\vcenter{\hbox{\input{figs/GD2VD2.5col-07.pstex_t}}} \to \ \vcenter{\hbox{\input{figs/GD2VD3col-07.pstex_t}}} $$
The next step meets obstruction, as we need to penetrate through arcs
that have been drawn. But let us continue neglecting the obstructions!
$$\vcenter{\hbox{\input{figs/GD2VD4col-07.pstex_t}}} \to \ \vcenter{\hbox{\input{figs/GD2VD5col-07.pstex_t}}}
\to \ \vcenter{\hbox{\input{figs/GD2VD6col-07.pstex_t}}} $$

What is obtained looks like a knot diagram, but, besides usual crossings,
it has double points which are not decorated. Such diagrams are called
{\sfit virtual link diagrams.\/} They were introduced by Kauffman in the
middle of nineties. Undecorated double points are called {\sfit virtual
crossings.\/}

In the construction above, virtual crossings emerged inevitably. The
only feasible way to avoid them is to attach handles to the plane and
use them as bridges.
$$\vcenter{\hbox{\input{figs/VD2kDonS-07.pstex_t}}} $$

\subsection{Link Diagrams on Orientable Surfaces}\label{s2.3}
A link diagram drawn on orientable surface $S$,
instead of the plane, defines a link in a thickened surface $S\times I$.
It defines a Gauss diagram, as well. 

{\sl Any Gauss diagram appears in this way. \/}  

This is proven by the construction above, with handles added when 
needed.\qed

{\sl For each Gauss diagram, there is the smallest orientable closed 
surface with a link diagram defining this Gauss diagram.}

Here by {\it smallest\/} I mean a surface with the greatest Euler
characteristic, but without components disjoint from the link diagram. 
To eliminate a cheap possibility of making the Euler
characteristic arbitrarily large by adding disjoint empty spheres, let
us require that each connected component of the surface contains a
piece of the link diagram.

To construct an orientable closed surface with greatest Euler
characteristic accommodating a link diagram with a given Gauss diagram,
one can first construct a germ of a link diagram on an oriented compact 
surface, which would contain the diagram as a deformation retract,
and then cap each boundary component of the surface with a disk. 
For details, see \cite{KK}.\qed 

A virtual link diagram may emerge as a projection to a plane of a link
diagram on an orientable surface embedded in $\R^3$.

\subsection{Moves}\label{s2.4}
\begin{figure}[hbt]
\centerline{\includegraphics
{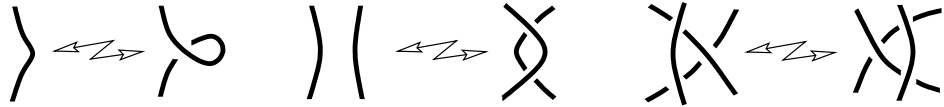}}
\caption{Reidemeister moves.}
\label{moves1}
\end{figure}

What happens to a link diagram when the link moves? It moves, too. 
A generic isotopy of a link can be decomposed into a sequence of isotopies 
each of which changes the diagram either as an isotopy of the plane 
or as one of three Reidemeister moves (see Figure \ref{moves1}). 
It does not matter, if the link lies in $\R^3$ and its diagram lies on $\R^2$,
or the link lies in a thickened oriented surface and its diagram lies on the
surface.

\begin{figure}[htb]
\centerline{\includegraphics
{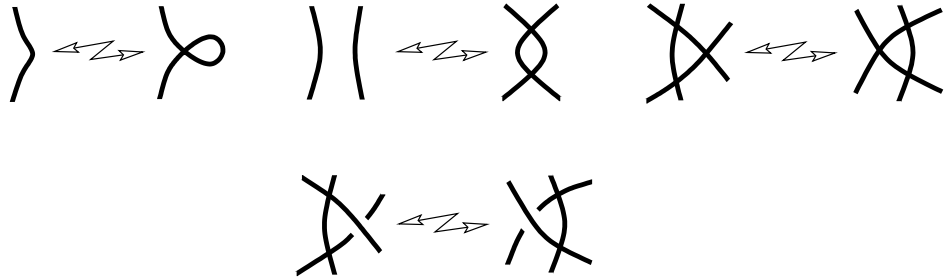}}
\caption{Virtual moves.}
\label{moves2}
\end{figure}

A virtual link diagram, which appeared as a plane projection 
of a link diagram on a surface, moves also as shown in Figure
\ref{moves2} when the link moves generically in the thickened surface.

All virtual moves can be replaced by detour moves:
\begin{figure}[hbt]
\centerline{\includegraphics[scale=.6]{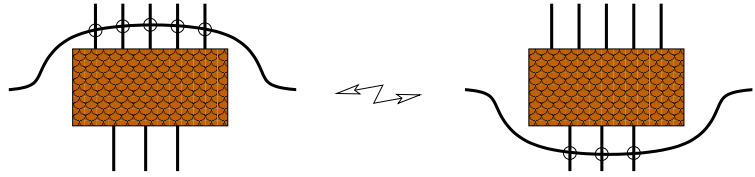}}
\caption{Detour move.}
\label{detour}
\end{figure}

Gauss diagram does not change under virtual moves.
Reidemeister moves act on Gauss diagram:  

\begin{table}[htbp]
\renewcommand{\arraystretch}{1.4}
\begin{tabular}{|m{2cm}|m{3.4cm}|m{6.8cm}|}
\hline Move's name & Reidemeister move &\hfil Its action on Gauss diagram\hfil \\ 
\hline
{Posi\-tive
first move}
&\includegraphics[scale=.8]{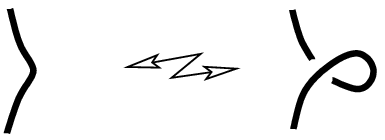}&\vspace{3pt}\hfil{\begin{picture}(0,0)%
\includegraphics{1+GDM.pstex}%
\end{picture}%
\setlength{\unitlength}{3947sp}%
\begingroup\makeatletter\ifx\SetFigFont\undefined%
\gdef\SetFigFont#1#2#3#4#5{%
  \reset@font\fontsize{#1}{#2pt}%
  \fontfamily{#3}\fontseries{#4}\fontshape{#5}%
  \selectfont}%
\fi\endgroup%
\begin{picture}(1997,646)(4478,-684)
\put(6226,-436){\makebox(0,0)[rb]{\smash{{\SetFigFont{12}{14.4}{\familydefault}{\mddefault}{\updefault}{\color[rgb]{0,0,0}$+$}%
}}}}
\end{picture}%
}\hfil\\
\hline
{Nega\-tive
first move}
&\includegraphics[scale=.8]{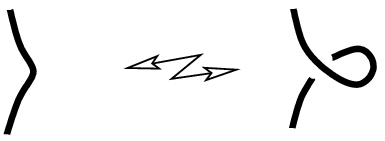}&\vspace{3pt}\hfil\begin{picture}(0,0)%
\includegraphics{1-GDM.pstex}%
\end{picture}%
\setlength{\unitlength}{3947sp}%
\begingroup\makeatletter\ifx\SetFigFont\undefined%
\gdef\SetFigFont#1#2#3#4#5{%
  \reset@font\fontsize{#1}{#2pt}%
  \fontfamily{#3}\fontseries{#4}\fontshape{#5}%
  \selectfont}%
\fi\endgroup%
\begin{picture}(2072,646)(4478,141)
\put(6301,389){\rotatebox{360.0}{\makebox(0,0)[rb]{\smash{{\SetFigFont{12}{14.4}{\familydefault}{\mddefault}{\updefault}{\color[rgb]{0,0,0}$-$}%
}}}}}
\end{picture}%
\hfil\\
\hline
{Second move}
&\vspace{3pt}{\includegraphics[scale=.8]{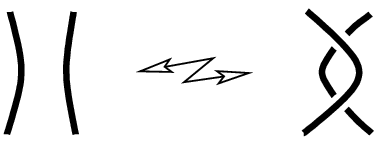}}&\vspace{3pt}{\input{figs/2GDM-70.pstex_t}}\\
\hline 
{Third move}
&{\includegraphics[scale=.8]{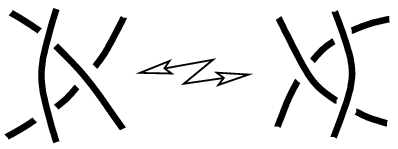}}&\vspace{8pt} {\input{figs/3GDM-70.pstex_t}}\\ 
\hline
\end{tabular}\vspace{10pt}
\caption{Action of Reidemeister moves on Gauss diagram.}
\label{t.1}
\end{table}

\subsection{Three Incarnations of Virtual Link Theory}\label{s2.5}
On the first stages of its development, the classical link theory 
received a completely combinatorial setup. Links were represented by
link diagrams and isotopies were represented by Reidemeister moves.

We see that virtual link theory has even {\bf two\/} similar 
{\bf combinatorial incarnations.\/} Virtual links (whatever they are) are represented, on one
hand, by {\it virtual link diagrams,\/} on the other hand, by 
{\it Gauss diagrams.\/} Furthermore, virtual isotopies (whatever they
are) are represented, on one hand, by {\it Reidemeister and detour
moves,\/} on the other hand, by the moves of Gauss diagrams
corresponding to Reidemeister moves.

{\bf Third incarnation,\/} truly topological one, is provided by
Greg Kuperberg \cite{Kup}. He proved that 

{\sl Virtual links up to virtual isotopy is the same as
irreducible links in thickened orientable surfaces up to orientation 
preserving homeomorphisms. \/}

Here {\sfit irreducible\/} means that each connected component of the 
thickened surface contains some part of the link and
it is impossible to find a simple curve $C$ on the surface, which is 
\begin{itemize}
\item disjoint from the projection of the link and
\item either non-zero homologous on the surface, or separating 
two parts of the link projection from each other, 
\end{itemize}
and it will be still impossible after any isotopy of the link.

In fact, Kuperberg \cite{Kup} described more specifically how a 
virtual link diagrams can be turned into an irreducible link in a 
thickened surface. He proved that a link diagram on an oriented 
surface can be {\sfit destabilized\/} to a link diagram of irreducible 
link on an oriented surface.  
A destabilization consists of embedded Reidemeister moves and 
Morse modifications of index 2 of the surface along a circle 
disjoint from the diagram. 
 
Kuperberg's results bridge combinatorics (=one-dimensional topology)
with the three-dimensional topology. The bridge can be used in both
directions: both for extending combinatorial techniques like quantum
link polynomials and link homology to links in 3-manifolds different
from $S^3$, and for using traditional topological techniques, like 
signatures, in the combinatorial environment.

For instance, an oriented link in a thickened surface realizes a 
homology class. A homeomorphism maps a link homologous to zero to a link
homologous to zero. Therefore the property of being homologous to zero
is a property of the virtual link. The same holds true for many other 
properties such as being homologous to zero modulo any number. 

Further, for a link homologous to zero modulo two in a thickened 
oriented surface one can define a {\it link signature.\/} 
Hence one can expect that there is a purely combinatorial 
construction of signature for virtual links of this kind.

\subsection{Twisted virtual links}\label{s2.6}
\begin{figure}[htb]
\centerline{$\fig{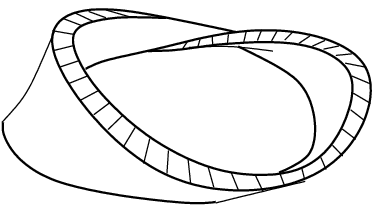}$}
\caption{}
\label{moebius-th}
\end{figure}

Non-orientable surface can be also thickened to an oriented 3-manifold.
For example, take  a M\"obius band $M$ embedded in $\R^3$ and thicken 
it, that is take its regular neighborhood. See Figure \ref{moebius-th}.
A neighborhood of $M$ in $\R^3$ is orientable and fibers over $M$.

To construct such a thickening, there is no need to embed 
a non-orientable surface to $\R^3$. Moreover, many of non-orientable
surfaces cannot be embedded into $\R^3$, but each of them has an 
orientable thickening. 
To thicken non-orientable surface $S$:
\begin{enumerate} 
\item Find an {\it orientation change curve\footnote{A curve realizing
homology class Poincar\'e dual to the first Stiefel-Whitney
class of $S$.}\/} $C$ 
(like {\it International date line\/}) on $S$.
$$\vcenter{\hbox{\begin{picture}(0,0)%
\hbox{\includegraphics{moebius.pstex}}%
\end{picture}%
\setlength{\unitlength}{2486sp}%
\begingroup\makeatletter\ifx\SetFigFont\undefined%
\gdef\SetFigFont#1#2#3#4#5{%
  \reset@font\fontsize{#1}{#2pt}%
  \fontfamily{#3}\fontseries{#4}\fontshape{#5}%
  \selectfont}%
\fi\endgroup%
\begin{picture}(2847,1738)(442,-1201)
\put(1306,-403){\makebox(0,0)[lb]{\smash{{\SetFigFont{7}{8.4}{\familydefault}{\mddefault}{\updefault}{\color[rgb]{1,0,0}$C$}%
}}}}
\end{picture}%
}}$$
\item  Cut $S$ along $C$: \ $S\mapsto S\fig{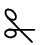}C$ \ 
$\fig{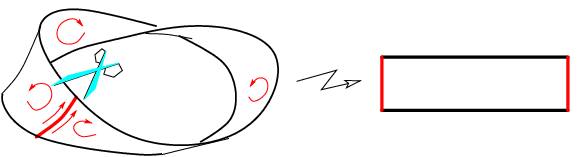}$
\item  Thicken: $(S\fig{figs/nozhnucw.eps}C)\times\R$.
\item Paste over the sides of the cut $(x_+,t)\sim (x_-,-t)$.
\end{enumerate}
Although the construction seems to depend on a choice of orientation
change curve, the thickening does not: it is unique up to a
homeomorphism and can be described as the total space of line bundle
over the surface with the first Stiefel - Whitney class equal to the
Stiefel - Whitney class of the surface.

A link in an orientable thickening of a non-orientable surface has a 
diagram on the surface. Since the fibration of the thickened surface
is not orientable one should take a special care on overpasses and 
underpasses. To distinguish them, one should orient the fiber over the
crossing. Since it is impossible to orient all fibers coherently, the
most natural solution is to use orientation of the fibration 
over the complement of an orientation change curve and keep the 
curve also shown on the diagram.

To encode also orientation of the thickening in the diagram, one has to
orient the complement of an orientation change curve on the base
surface. The local orientation of the base and the global orientation 
of the total space determine an orientation of a fiber. It is this
orientation, which is used to distinguish overpasses and underpasses.

An orientation of the complement of an orientation change curve on the base
surface is not determined by the link in oriented thickening of a
non-oriented surface. It can be reversed on any component of the
surface. The reversing switches over- and under-passings on this
component, but preserves local writhes which are defined as follows.
$$\vcenter{\hbox{\input{figs/LocWr.pstex_t}}}$$

A generic isotopy of a link in an orientable thickening of a
non-orientable surface decomposes into a sequence of isotopies each of
which acts on the diagram either as an isotopy of the surface, or 
a Reidemeister move, or one of the following two additional moves, in
which the orientation change curve is involved:

\centerline{$\fig{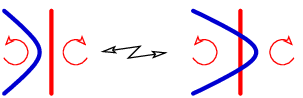}$ \ and \ $\fig{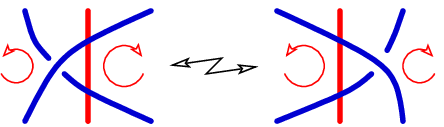}$}

The first of these moves happens when a piece of link penetrates 
through the the preimage of the orientation change curve. The second
happens when a crossing moves through the orientation change curve.
Since the orientation of of the fiber changes at the moment, the
overpass and underpass exchange. 

In practice, it is convenient to cut the diagram along the orientation
change curve, but keep in mind the identification which would allow to
recover the cut. Say, in the case of thickened projective plane, for
an orientation change curve can be taken a projective line. The cut
along it gives a disk, which is much more comfortable to draw on, than
the projective plane. In this case, the theory sketched above turns 
into the theory of diagrams for links in the projective space developed
by Julia Drobotukhina \cite{D1}. 

In general case links in an oriented thickening of a non-orientable
surface were considered by Mario O.Bourgoin in his Research Statement 
\cite{B}. He suggested also the corresponding generalization of virtual
links, twisted virtual links, and announced a generalization of
Kuperberg's theorem to the twisted setup.

We consider the corresponding generalization of Gauss diagrams (rather
than generalization of virtual link diagrams outlined by Bourgoin
\cite{B}). 

{\it Twisted Gauss diagram\/} is a Gauss diagram with a finite set of
points marked on the base curve. No marked point coincide with an
end-point of an arrow. 

The two extra moves of diagrams considered above correspond to the
following moves:

\centerline{$\fig{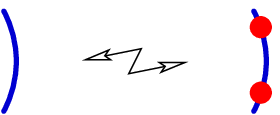}$ \ \qquad and \ \qquad $\vcenter{\hbox{\input{figs/DrobMove.pstex_t}}}$}

\subsection{Stripping of the Third Dimension}\label{s2.7}
Recall that a classical link diagram is a decorated generic projection 
of the link. The generic projection is a generic immersion. A great
piece of what was said above about link diagrams can be repeated, with
appropriate simplifications, about generic immersions of 1-manifold. For
short, we call the image of a generic immersion of a closed 1-manifold 
to a surface $S$ a {\it generic curve\/} on $S$. The immersion gives
rise to a parametrization of the curve.   

A closed 1-manifold with a finite number of chords connecting pairwise
distinct points of the 1-manifold is called a {\sfit chord diagram.\/}
The 1-manifold is called the {\sfit base\/} of the chord diagram.
A chord diagram, each chord of which is equipped with orientation, is
called an {\sfit arrow diagram.\/}
A chord diagram, each chord of which is equipped with a sign, is called
a {\sfit signed chord diagram.\/}
Gauss diagrams considered above are signed arrow diagrams.

A generic immersion of a 1-manifold to a surface defines chord diagram,
in which the base is the source 1-manifold and each chord connects
points having the same image. If the source 1-manifold and target 
surface of the immersion are oriented, the chords get natural orientations:
direct a chord from branch $A$ to branch $B$ such that the orientation
at the target is defined by the basis formed of vectors which are the 
images of tangent vectors to $A$ and $B$ defining the orientation of
the source and taken in this order.   

In the case of a link diagram, these orientations of chords can be
obtained from the orientations and local writhes involved in the Gauss 
diagram by multiplying them: if the sign of the arrow in Gauss diagram is
$+$, take its orientation intact, if the sign is $-$, reverse the
orientation. 

Not any arrow diagram and even chord diagram can be generated by a
generic immersion to plane. It was the problem attracted Gauss how to
recognize which diagrams appear in this way. This problem received 
solutions in a number of ways, but we refrain from going into this 
vast matter.

A step from generic curves on an oriented surface parallel to the step 
from link diagrams on an oriented surface 
to virtual link diagrams gives rise to {\sfit flat virtual knots\/}, 
see D.Hrencecin, L. Kauffman \cite{HK}.
The counter-part of Reidemeister moves are {\sfit flat Reidemeister
moves\/}, see \cite{HK}.

I am not aware about any counter-part of the Kuperberg Theorem, which
would relate flat virtual knots considered up to flat Reidemeister moves 
to irreducible generic immersions of a 1-manifold to an orientable surface 
considered up to homeomorphism and homotopy.   

\section{Kauffman Bracket of Virtual Links}\label{s3}

\subsection{Digression on Kauffman Bracket of Classical Link}\label{s3.1}
Kauffman bracket of a classical link diagram is a Laurent polynomial 
in $A$ with integer coefficients
     $$\langle\text{Link diagram}\rangle \in \Z[A,A^{-1}]$$

For example,
$$
\begin{aligned}
&\langle\text{unknot}\rangle= \ 
       &{\langle\bigcirc\rangle= \ }
       &{-A^2-A^{-2}}\\
&{\langle\text{Hopf link}\rangle= \ }
       &{\langle\vcenter{\hbox{
	\includegraphics[scale=.2]{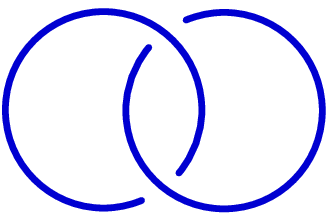}}} \ \rangle= \ }
       &{A^6+A^2+A^{-2}+A^{-6}}\\
\end{aligned}
$$
$$
\begin{aligned}
&{\langle\text{empty link}\rangle= \ }
       &{\langle \ \ \rangle= \ }
       &{1}\\
&{\langle\text{trefoil}\rangle= \ }
       &{\langle\vcenter{\hbox{
	\includegraphics[scale=.4]{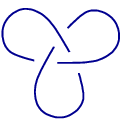}}} \ \rangle= \ }
       &{A^7+A^3+A^{-1}-A^{-9}}\\
&{\langle\text{figure-eight knot}\rangle= \ }
       &{\langle\vcenter{\hbox{
	\includegraphics[scale=.2]{figs/d-fig8.eps}}} \ \rangle= \ }
       &{-A^{10}-A^{-10}}
\end{aligned}
$$

Kauffman bracket is defined by the following properties:
\begin{enumerate}
\item {  
$\langle\bigcirc\rangle=-A^2-A^{-2}$,} 
\item { $ \langle D\ \amalg \ \bigcirc \rangle 
=(-A^2-A^{-2})\langle D\rangle$,} (here $\amalg$ means disjoint union)
\item {
$\langle\vcenter{\hbox{\includegraphics{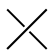}}}\rangle=
A\langle\vcenter{\hbox{\includegraphics{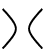}}}\rangle +
A^{-1}\langle\vcenter{\hbox{\includegraphics{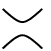}}}\rangle $
({\it Kauffman Skein Relation\/})
.} 
\end{enumerate}
                                    
Indeed, applying the last property to each crossing of a link diagram, 
one reduces the diagram to collections of embedded circles.
Then the first two properties complete the job.
This calculation can be summarized in the following 
{\sfit state Sum Model.}

A {\sfit  state\/} of diagram is a distribution of {\sfit markers\/} 
over all crossings.  At each crossing of the
diagram there should be a marker specifying a pair of vertical angles:
$\vcenter{\hbox{\includegraphics{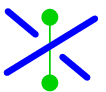}}}$, \ 
 $\vcenter{\hbox{\includegraphics{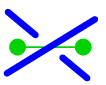}}}$.

For example, knot diagram:
$\vcenter{\hbox{\includegraphics[scale=.5]{figs/d-fig8.eps}}}$ has 
states: 
{$\vcenter{\hbox{\includegraphics[scale=.5]{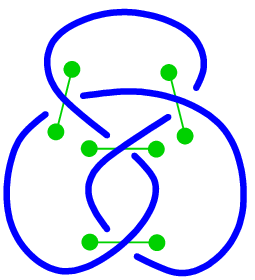}}}$\ ,}
{ $\vcenter{\hbox{\includegraphics[scale=.5]{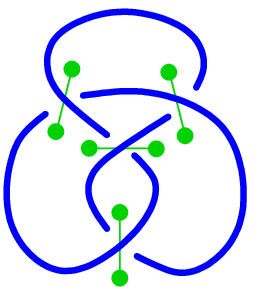}}}$\ ,}
{ $\vcenter{\hbox{\includegraphics[scale=.5]{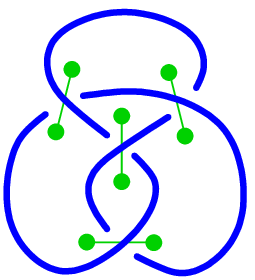}}}$,}
{$\vcenter{\hbox{\includegraphics[scale=.5]{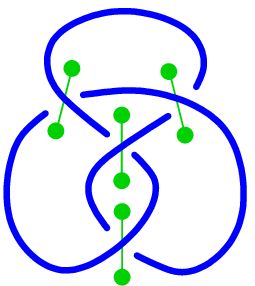}}}$\ ,}
{\dots} 
{Totally $2^c$ states, where $c$ is the number of crossings.} 
   
Three numbers are associated to a state $s$ of diagram $D$: \begin{itemize}
\item the number $a(s)$ of {\sfit positive\/} markers 
$\vcenter{\hbox{\includegraphics{figs/marker+.eps}}}$,
\item the number $b(s)$ of {\sfit negative\/} markers 
$\vcenter{\hbox{\includegraphics{figs/marker-.eps}}}$,
\item {the number $|s|$ of components of the curve $D_s$ obtained 
by smoothing of $D$ along the markers of $s$.}
\end{itemize}
For example, state
{ $s = \ \vcenter{\hbox{\includegraphics[scale=.5]{figs/d-fig8.m1.eps}}}
$} has $a(s)=1$, \ $b(s)=3$, { $\text{smoothing}(s)=
\vcenter{\hbox{\includegraphics[scale=.5]{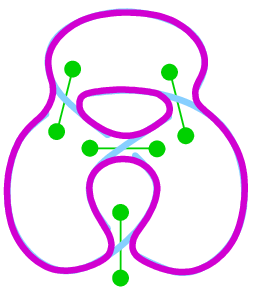}}}$}, 
and { \ $|s|=2$}.

The contribution of a state $s$ to Kauffman bracket along with the
calculation sketched above is $A^{a(s)-b(s)}(-A^2-A^{-2})^{|s|}$.
Finally, the whole Kauffman bracket is equal to the following
{\sfit:state sum\/}  
$$\langle D\rangle= \sum_{s\text{ state of }D}
A^{a(s)-b(s)}(-A^2-A^{-2})^{|s|}$$                             

Example: consider Hopf link, 
$\vcenter{\hbox{
	\includegraphics[scale=.4]{figs/HopfL.eps}}}$\\
 {$\left\langle\vcenter{\hbox{
	\includegraphics[scale=.4]{figs/HopfL.eps}}} \ \right\rangle= \ $\\ 
       \vskip3pt }
 {$\left\langle\vcenter{\hbox{
	\includegraphics[scale=.4]{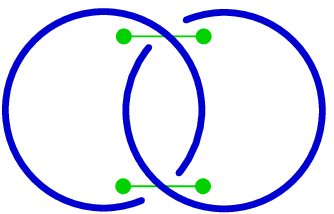}}} \ \right\rangle +
	\left\langle\vcenter{\hbox{
	\includegraphics[scale=.4]{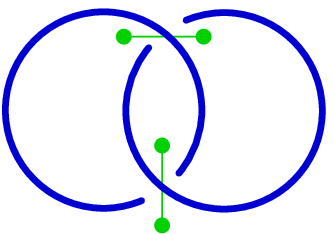}}} \ \right\rangle +
	\left\langle\vcenter{\hbox{
	\includegraphics[scale=.4]{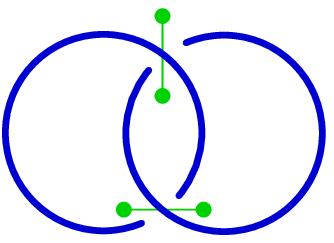}}} \ \right\rangle +
	\left\langle\vcenter{\hbox{
	\includegraphics[scale=.4]{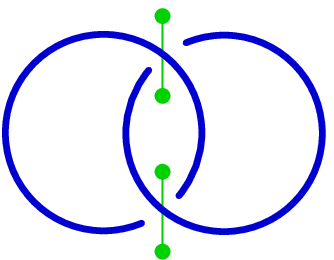}}} \ \right\rangle = \ $\\ 
	\vskip3pt } 
\noindent
{$A^2(-A^2-A^{-2})^2 + 2(-A^2-A^{-2})  +
 A^{-2}(-A^2-A^{-2})^2 = $ \\ }
{$(A^6 +2A^2+A^{-2}) -2A^2 -2A^{-2}+(A^2+2A^{-2}+A^{-6}) = $ \\ }
{$A^6 + A^2 + A^{-2} +A^{-6} $. }
                             
\subsection{Kauffman State Sum in Terms of Gauss Diagram}\label{s3.2}

Let us rewrite the state sum above in terms of Gauss diagram.
Recall that to a crossing of a classical link diagram there corresponds
an arrow equipped with a sign in the corresponding Gauss diagram.
$$\vcenter{\hbox{\input{figs/cr-ar.pstex_t}}}$$ 

Clearly, to a smoothing of a crossing there corresponds a surgery along
the corresponding arrow:
\\  \vspace{8pt}
$$\vcenter{\hbox{\input{figs/cr-ar.sm+.pstex_t}}}$$  
$$\vcenter{\hbox{\input{figs/cr-ar.sm-.pstex_t}}}$$  

We see that the type of smoothing depends only on two signs, the sign of
the marker and the sign of the crossing. Namely, the surgery preserves 
orientation, if these signs coincide; if the signs are opposite,
the surgery reverses orientation. 

Thus, a state and all the numerical characteristic of a state involved
in the Kauffman state sum can be read from the
Gauss diagram. A state of Gauss diagram is a distribution of signs 
(signs of markers) over the set of all of its arrows. The number $|s|$
of components of the curve obtained by the smoothing along markers
can be figured out from the Gauss diagram and signs of markers.

Moreover, for recovering of the Kauffman state sum we even do not need
an important ingredient of Gauss diagram, directions of its arrows.
Therefore we can forget about them. Let us check if this is a useful 
possibility.

\subsection{Blunted Gauss Diagrams}\label{s3.3}
A Gauss diagram with directions of arrows forgotten is a signed chord
diagram. The forgetting of directions is called {\sfit blunting\/}, 
the result is called a {\sfit blunted Gauss diagram.\/} 

Moves of Gauss diagrams defined by Reidemeister moves under blunting
turns to the following moves of signed chord diagrams:

$$\vcenter{\hbox{\begin{picture}(0,0)%
\includegraphics{1BGDM.pstex}%
\end{picture}%
\setlength{\unitlength}{2368sp}%
\begingroup\makeatletter\ifx\SetFigFont\undefined%
\gdef\SetFigFont#1#2#3#4#5{%
  \reset@font\fontsize{#1}{#2pt}%
  \fontfamily{#3}\fontseries{#4}\fontshape{#5}%
  \selectfont}%
\fi\endgroup%
\begin{picture}(1997,646)(4478,-684)
\put(6301,-436){\makebox(0,0)[rb]{\smash{{\SetFigFont{7}{8.4}{\familydefault}{\mddefault}{\updefault}{\color[rgb]{0,0,0}$\pm$}%
}}}}
\end{picture}%
}},\qquad\qquad \vcenter{\hbox{\input{figs/2BGDM.pstex_t}}},$$
\vspace{5pt}
$$\vcenter{\hbox{\input{figs/3BGDM.pstex_t}}}$$   
We shall call these moves {\sfit Reidemeister moves of signed chord
diagrams.\/}

As shown above, Kauffman state sum is defined in terms of the
corresponding blunted Gauss diagram. Such a state sum is defined 
for any signed chord diagram.
The classical proof of invariance of the Kauffman bracket under 
the second and third Reidemeister moves works perfectly in the setup
of signed chord diagrams. Moreover, under the first Reidemeister move
the Kauffman bracket of a signed diagram behaves exactly as in the
classical setup: the positive first move causes multiplication of 
the Kauffman bracket by $-A^3$, the negative one, by $-A^{-3}$.

Thus, the Laurent polynomial 
$$f_D(A)=(-A)^{-3w(D)}\langle D \rangle,$$ 
where $w(D)$ is the sum of signs of all chords of
a signed chord diagram $D$, is invariant under all Reidemeister moves. 
As well-known, it is closely
related to the Jones polynomial in the classical case. However, as
everything here works just fine for blunted Gauss diagrams, it can be 
(and it was) taken as a definition of the Jones polynomial for virtual 
links.

\subsection{Twisted Versus Blunted}\label{s3.4}
Observe that under each of the additional two moves of twisted Gauss 
diagrams the signs of arrow do not change. 

Therefore if we forget in a twisted Gauss diagram both directions of
arrows and points, we get a signed chord diagram and moves of twisted Gauss
diagrams turn into Reidemeister moves of signed chord diagrams.

This gives Kauffman bracket and Jones polynomial for twisted virtual
links. This has been done in literature, as well. For links in oriented 
thickened projective plane (which is equivalent to the 3-dimensional 
projective space $\R P^3$, as the oriented thickened projective plane 
is the complement of a point in $\R P^3$), 
the Kauffman bracket and Jones polynomial were defined and studied
by Drobotukhina \cite{D1} in 1990. For general twisted virtual links
these polynomials and their refinement was outlined by Bourgoin
\cite{B}.

At first glance, everything is preserved when we pass from the
classical links to virtual and even twisted virtual ones.
However, this impression is a little bit misleading. 
Some properties of the Kauffman bracket and Jones polynomial 
change drastically.

\subsection{Exponents in Kauffman Bracket}\label{s3.5}
As we have seen in the examples
  $$
\begin{aligned}
&\langle\text{unknot}\rangle&= \ 
       &\langle\bigcirc\rangle&= \ 
       &-A^2-A^{-2}\\
&\langle\text{Hopf link}\rangle&= \ 
       &\langle\vcenter{\hbox{
	\includegraphics[scale=.2]{figs/HopfL.eps}}} \ \rangle&= \ 
       &A^6+A^2+A^{-2}+A^{-6}\\
&\langle\text{empty link}\rangle&= \ 
       &\langle \ \ \rangle&= \ 
       &1\\
&\langle\text{trefoil}\rangle&= \ 
       &\langle\vcenter{\hbox{
	\includegraphics[scale=.4]{figs/trefoil.eps}}} \ \rangle&= \ 
       &A^7+A^3+A^{-1}-A^{-9}\\
&\langle\text{figure-eight knot}\rangle&= \ 
       &\langle\vcenter{\hbox{
	\includegraphics[scale=.2]{figs/d-fig8.eps}}} \ \rangle&= \ 
       &-A^{10}-A^{-10},
\end{aligned}
$$ 
exponents in Kauffman bracket of a classical link are congruent to
each other modulo 4.  

 This happens because: \begin{itemize} 
\item the contribution 
$A^{a(s)-b(s)}(-A^2-A^{-2})^{|s|}$ of each state has this property,
\item from each state we can get to any other state changing one marker
a time,
\item  and changing a single marker in a state changes
$a(s)-b(s)$ by 2 and $|s|$ by 1. 
\end{itemize}

However, this is not the case for virtual links. For example,
$$\left\langle\fig{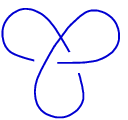}\right\rangle \ = \
A^{-4}+A^{-6}-A^{-10}$$

$$\left\langle\fig{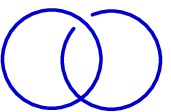}\right\rangle \ = \ A+A^{-1}$$  

What is wrong in the proof above if the link is virtual? 
\begin{itemize} 
\item  Still, in the contribution  $A^{a(s)-b(s)}(-A^2-A^{-2})^{|s|}$ of 
each state all the exponents are congruent to each other modulo 4. 
\item  Still, from each state we can get to any other state changing 
one marker a time.
\item   Still, changing a single marker changes $a(s)-b(s)$ by two.
However, $|s|$ may be preserved.  
\end{itemize}

Change of a single marker causes a Morse modification of the result 
of smoothing. 
In the classical case this {\it Morse modification is embedded in plane,
and therefore preserves orientation.\/} 

As well-known, {\it a preserving orientation 
Morse modification of a 1-manifold  changes the number of components
by one. A Morse modification, which does not preserve orientation
of a connected 1-manifold, preserves the number of connected components.\/}
Morse modification non-preserving orientation, would cause a shift of 
exponents by 2.

\section{Orientations of Chord Diagrams}\label{s4}

\subsection{Orientation of Chord Diagram}\label{s4.1}
What structure on blunted Gauss diagram would guarantee that, for any
of its states, the corresponding smoothing have an orientation
such that change of any marker causes a Morse modification of the result
of smoothing {\it preserving\/} the orientation? 

\begin{figure}[htb]
\centerline{\epsffile{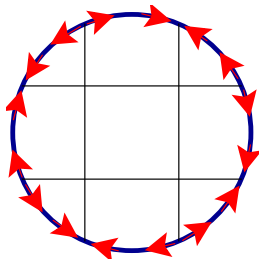}}
\caption{Orientation of a chord diagram.}
\label{AltOr}
\end{figure}

Such a structure is a collection of orientations of arcs of the 
base 1-manifold between end points of chords such that these
orientations define an orientation of each smoothing. 
For this the orientations should have the
following two properties:
\begin{itemize}
\item The orientations cannot be extended over an end point of a
chord,
\item For each chord, one of its end points is attractive, while the
other one is repulsive. 
\end{itemize}
See Figure \ref{AltOr}.

The first property means that the orientations of arcs  
induce anonymously an orientation of each end point of each chord. 
The second property means that for each chord the induced orientations 
of its end points are opposite to each other.

Orientations of arcs of a chord diagram satisfying these two properties
is called an {\sfit orientation\/} of chord diagram.

\begin{figure}[htb]
\centerline{\epsffile{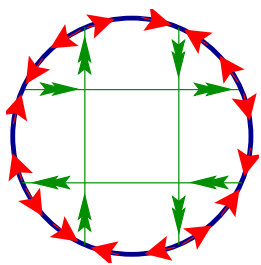}}
\caption{Orientations of chords in an oriented chord diagram.}
\label{AltOrCh}
\end{figure}

Orient chords so that the induced orientations of their end points were
opposite to the orientations induced by orientations of adjacent arcs.
See Figure \ref{AltOrCh}.
 
Chords become arrows. Orientations of the arrows alternate in the sense
that along the base of the diagram each arrowhead is followed by an 
arrowtail and each arrowtail is followed by arrowhead. 

Observe that an orientation of a chord diagram defines an orientation 
on the result of any of its smoothings. Indeed,
the result of a smoothing is composed from arcs of the base and arcs 
which go along chords. The orientations of these pieces agree with each
other.

\subsection{Alternatable Virtual Link Diagrams}\label{s4.2}
Observe that directions of arrows in a Gauss diagram of a classical alternating 
link diagram alternate in the same way as the orientations of chords
considered above.

Recall that a classical link diagram is called {\sfit
alternating,\/} if along a branch of the link over-crossing would 
always follow after under-crossing and under-crossing would always 
follow after over-crossing.

Similarly, a virtual link is called {\sfit alternating\/}, if along 
a branch of the link over-crossing always follow after under-crossing 
and under-crossing always follow after over-crossing.
By switching some overpasses and underpasses, a virtual link diagram 
such that its blunted Gauss diagram is orientable, can be made
alternating. Vice versa, if, by switching some overpasses and underpasses, 
a virtual link diagram can be made alternating, its blunted Gauss
diagram is orientable.

Therefore a virtual link diagram with orientable blunted Gauss diagram is 
called {\sfit alternatable\/}.

\subsection{Moves of Chord Diagrams}\label{s4.2.5}
\begin{figure}
$$\sfig{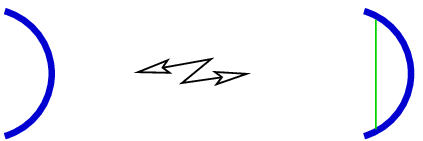}{.6},\qquad\qquad \sfig{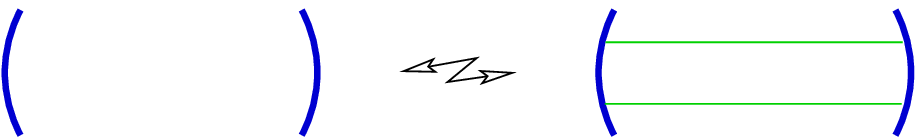}{.6},$$
\vspace{5pt}
$$\sfig{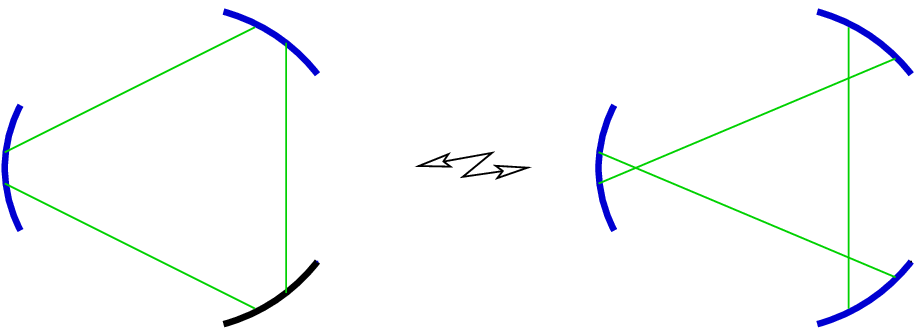}{.6}$$   
\caption{Reidemeister moves of chord diagrams}
\label{RMChD}
\end{figure}

Transformations of chord diagram shown in Figure \ref{RMChD} are called
{\sfit Reidemeister moves\/}.

\begin{thm}\label{4.7.A}
The result of any first or third Reidemeister move applied to an 
orientable chord diagram is orientable. 
\end{thm}

\begin{figure}[htb]
\centerline{$\fig{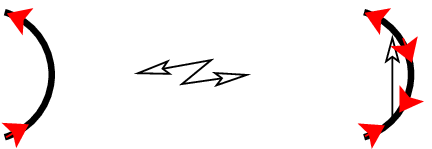}$}
\centerline{First Reidemeister move.}\vskip.1cm
\centerline{$\fig{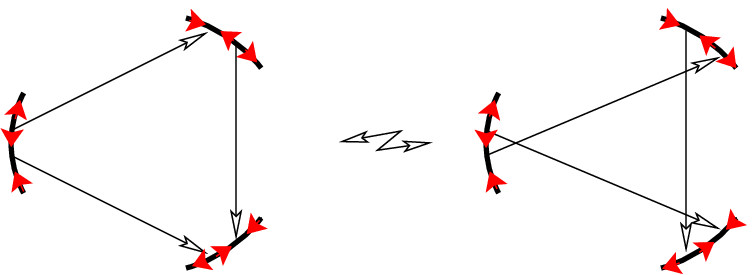}$}
\centerline{Third Reidemeister move.}
\caption{Behavior of orientation of a signed chord diagram under
first and third Reidemeister moves.}
\label{OrOfGDM}
\end{figure}

\begin{proof} An orientation of a chord diagram prior to the
move in the fragment which is about to change is unique up to reversing.
Its replacement admits an orientation coinciding with the
original one near the boundary. See Figure \ref{OrOfGDM}.
\end{proof}

The second Reidemeister move can destroy orientability.
Of course, a second Reidemeister move can transform a orientable chord
diagram to an orientable one.  All
second Reidemeister moves decreasing the number of chords and 
about half of all second Reidemeister moves increasing it 
which can be applied to a chord diagram preserve orientability.  
One can easily recognize, if a second Reidemeister move preserves 
orientability, see Figure \ref{2OrGDM}.

\begin{figure}[htb]
\centerline{\input{figs/2OrChDM.pstex_t}}
\caption{Behavior of orientations of a chord diagram under
second Reidemeister moves.}
\label{2OrGDM}
\end{figure}

\subsection{Checkerboard Coloring of a Classical Link
Diagram}\label{s4.3}
 
Let us analyze where the orientation of a chord diagram underlying a
Gauss diagram of a classical link comes from. 

A generic curve on a plane can be considered as 
a 1-cycle with coefficients in $\Z_2$ on the plane. As a plane has 
trivial homology, this 1-cycle bounds a 2-chain. The latter is described 
as the union of all black domains in the checkerboard coloring of the 
diagram. 

\begin{figure}[h]
\centerline{\includegraphics[scale=1.5]{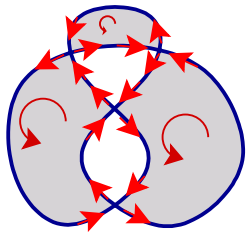}}
\caption{Checkerboard coloring and orientations.}
\label{f1}
\end{figure}

Each connected component of the complement of a generic curve inherits 
an orientation from the whole plane. These orientations of the black 
domains induce orientations on their boundaries. The boundary of a 
black domain consists of pieces of the curve. 
Thus an arc of the curve between any two consecutive double 
points gets a natural orientation. See
Figure \ref{f1}. It is called a {\sfit checkerboard orientation.\/}

Obviously, the checkerboard orientation gives orientation of the
corresponding chord diagram. The orientation
induced by the orientation of chord diagram on the result of a 
smoothing of the curve can be also obtained as a
checkerboard orientation. 

\begin{figure}[hbt]
\centerline{\includegraphics[scale=.9]{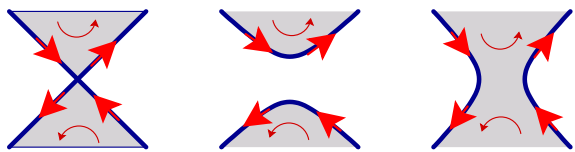}}
\caption{Checkerboard orientation of smoothings.}
\label{OrOfSmooth}
\end{figure}

Indeed, the result of a smoothing also admits 
a checkerboard coloring which coincides with the checkerboard coloring 
of the diagram outside of small neighborhoods of its double points. 
The orientation induced on the smoothing by the orientations of the 
black domains agree with the checkerboard orientation of the diagram. 
See Figure \ref{OrOfSmooth}. 

\subsection{On Orientable Surfaces}\label{s4.4}
The same arguments work for a link diagram zero-homologous modulo two 
on an oriented surface. An orientable surface with a zero-homologous 
generically immersed curve admits a checkerboard coloring. We see that:

\begin{thm}\label{T1}
The chord diagram of a generic curve on an orientable surface 
admitting a checkerboard coloring is orientable. 
\end{thm}

\begin{cor}[Naoko Kamada \cite{NK}] 
A link diagram on an orientable surface admitting a 
checkerboard coloring is alternatable. \qed  
\end{cor}

\begin{cor}[Naoko Kamada \cite{NK}]
Exponents of monomials in the Kauffman bracket of a 
checkerboard colorable link diagram on an oriented surface are congruent
to each other modulo 4.
\end{cor}\qed

It may happen that an alternatable link diagram on an orientable surface
does not admit a checkerboard coloring. Indeed, any diagram admitting a
checkerboard coloring can be spoiled by a Morse modification of index 1,
i.e., removing two disjoint open disks and attaching a tube connecting
there boundary circles. For this, take the disks in the complementary domains
colored with different colors. Of course, this stabilization of virtual
link diagram does not destroy alternatability, which is a property of
Gauss diagram. 

Nonetheless, according to the following theorem, which is also 
basically  due to Naoko Kamada \cite{NK}, this cannot happen to 
irreducible diagram.

\begin{thm}\label{T1^-1}
A generic curve on an orientable surface such that each connected 
component of its complement is a disk and its chord diagram is
orientable admits a checkerboard coloring. 
\end{thm}

\begin{cor}[Naoko Kamada \cite{NK}] 
An alternatable link diagram on an orientable surface
such that each connected component of its complement is a disk admits a
checkerboard coloring.
\end{cor} 

\begin{proof}[Proof of \ref{T1^-1}]
An orientation of the chord diagram gives rise to an orientation
of the boundary of each connected component of the complement of the
curve. Since each such component is a disk, the orientation of its
boundary orients the component itself. On the other hand, the
orientation of the whole surface induces an orientation on each of these
components. For some of the components these two orientations coincide,
for the others they are opposite to each other. The components for which
the orientations coincide form a chain modulo 2 bounded by the curve. 
By coloring these components in black and the others in
white, we get a desired checkerboard coloring of the surface. 
\end{proof}

The condition about components of the complement can be weakened:
instead of requiring that they are homeomorphic to disk, it suffices to
require that intersection of the curve with the closure of each of the
components is connected.

\subsection{Alternatable Virtual Links}\label{s4.7}
Recall that a virtual link diagram which gives rise to orientable
blunted Gauss diagram is called {\sfit alternatable.\/} 

\begin{thm}[Corollary of Theorem \ref{4.7.A}]\label{4.7.A-c}
The result of any first or third Reidemeister move applied to an 
alternatable virtual link diagram is alternatable. \qed
\end{thm}

\begin{figure}[htb]
\centerline{\epsffile{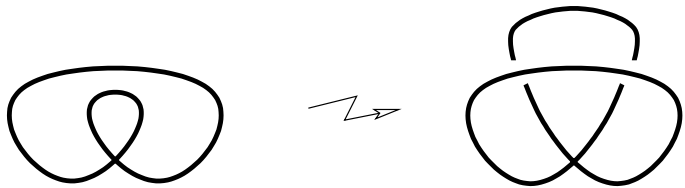}}
\caption{Creating a non-alternatable virtual knot diagram from a virtual 
diagram of unknot with single virtual crossing by a second 
Reidemeister move.}
\label{A2NA}
\end{figure}

A second Reidemeister move can turn an alternatable virtual link diagram 
to a non-alternatable, see Figure \ref{A2NA}.

Moreover, each alternatable virtual diagram can be made non-alternatable
by a single second Reidemeister move. Therefore alternatability is not 
a property shared by all virtual diagrams of a given virtual link.

Alternatable isotopy is a sequence of Reidemeister moves preserving 
alternatability and virtual moves of an alternatable virtual link
diagram. 

\begin{thm}\label{4.7.B}
Alternatable virtual link diagrams are virtualy isotopic, iff they 
can be related by an alternatable isotopy.
\end{thm}

\begin{proof} Let $D_1$ and $D_2$ be virtually isotopic alternatable
link diagrams. Realize each of them as a diagram admitting a 
checkerboard coloring on an oriented closed surface. 

According to the Kuperberg Theorem \cite{Kup}, virtually isotopic link 
diagrams can be destabilized to link diagrams on an oriented surface 
$S$, where they can be related by embedded moves. A destabilization 
consists of embedded Reidemeister moves and Morse modifications of 
index 2 of the ambient surface along a circle disjoint with the 
diagram. 

Neither an embedded Reidemeister move, nor Morse modification of index 
two  can destroy a checkerboard coloring. A Morse modification of index 
two does not change the Gauss diagram. Therefore the sequence of moves 
connecting the Gauss diagrams corresponding to $D_1$ and $D_2$ 
which correspond to the Reidemeister moves existing by the 
Kuperberg Theorem constitute an alternatable isotopy.
\end{proof}

Thus alternatable virtual links do not form a new category of objects
similar to virtual links, but are virtual links of special kind. 
They have special properties. For example, the exponents of Kauffman
bracket are congruent to each other modulo 4, see the second Corollary
of Theorem \ref{T1} above. From the point of view of 3-dimensional
topology, they can be characterized as irreducible links in 
thickened oriented surfaces which realize trivial homology class modulo
2.

\subsection{On a Non-Orientable Surface}\label{s4.5}
Theorem \ref{T1} and its corollary cannot be extended literally to 
twisted link diagrams on non-orientable surfaces . 
For example, link $2_1$ on the projective plane
has checkerboard colorable diagram, but its Kauffman bracket is
$A^{-4}+A^{-6}-A^{-10}$, see \cite{D2}.
\begin{figure}[htb]
\centerline{\includegraphics[scale=.5]{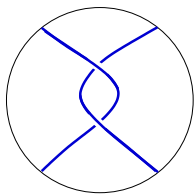}}
\caption{Projective link $2_1$}
\label{2_1}
\end{figure}

However, a generalization recovers, if the property of being
zero-homologous modulo two is generalized in other way.

\begin{thm}\label{T2}
The chord diagram of a generic curve realizing 
the homology class dual to the first Stiefel-Whitney class
of the surface is orientable.  
\end{thm}
 
\begin{proof} Let $C$ be a generic curve on a surface $S$ and
realizing a homology class dual to $w_1(S)\in H^1(S;\Z_2)$. The
complement $S\sminus C$ admits an orientation which cannot be extended 
over $C$ at any point. Each arc of $C$ is involved in $C$ with
multiplicity two in the boundary of the corresponding fundamental
class $[S\sminus C]$. At a double point the orientations of $S\sminus C$
and arcs of $C$ look as follows: 
$$\vcenter{\hbox{\includegraphics[scale=.5]{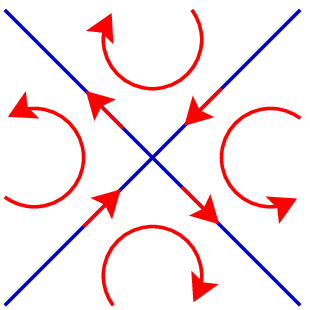}}}.$$
Thus the orientations of arcs of $C$ gives rise to an orientation of
the corresponding chord diagram.
\end{proof}

\begin{cor}Exponents of monomials in the Kauffman bracket of a link 
non-zero-homologous in the real projective space are congruent to each
other modulo 4.  \qed
\end{cor}

\begin{thm}[Generalization of Theorem \ref{T1^-1}]\label{T2^-1}
If the chord diagram of a generic curve on a surface is orientable
and each connected component of its complement on the surface is a disk 
then the curve realizes the homology class dual to the first 
Stiefel-Whitney class of the surface.  
\end{thm}

\begin{proof} An orientation of the chord diagram gives rise to an 
orientation of the boundary of each connected component of the 
complement of the curve. Since each such component is a disk, the 
orientation of its boundary orients the component itself. 

The orientations of the components of complement form an orientation 
of whole complement, which does not extend over the curve, because
from both sides of the curve it induces the same orientation on each
piece of the curve. Therefore the curve realizes the obstruction to
orientability of the surface.
\end{proof}

\subsection{Obstruction to Orientability of Chord Diagram}\label{s6.1}
Near each chord of a chord diagram an orientation of the diagram 
looks standard, up to simultaneous reversing.
For any choice of the orientation, on each arc of the base either 
the orientations at the end points extend to the whole arc, or not. 
Thus for any choice of orientations of chords there is a well 
defined $\Z/2\Z$-valued function of the set of arcs of the base.
It can be considered as a 1-cochain with values in $\Z/2\Z$ on the 
chord diagram.
Reversing of a single chord's orientation causes change of values on 
each of the four adjacent arcs. 

This means that there is a well-defined one-dimensional cohomology class
with values in $\Z/2\Z$ of the underlying space of the chord diagram
(the union of the base and chords).
This class is zero, iff the diagram is orientable. Therefore we call
it the {\sfit obstruction to orientability. \/}

The class is defined by topology of the pair consisting of the 
underlying space and the union of all chords. It is not defined 
by the homotopy type of the underlying space, although it belongs to a 
cohomology group depending only on the homotopy type. 
It is not defined by the topology of the four-valent 
graph obtained by contracting each of the chords, either. Indeed, 
it depends on division of arcs adjacent to a vertex to pairs 
of opposite arcs, that is arcs adjacent to one end point of 
the chord in the underlying space of the chord diagram. 

The simplest example is a circle with two chords intersecting each 
other, on one hand, and two circles connected with two chords, on the
other hand. The quotient four-valent graphs are homeomorphic.

\section{Khovanov Homology of Oriented Signed Chord Diagrams}\label{s5}

\subsection{Khovanov Homology of Classical Links}\label{s5.1}
In 1998 Khovanov \cite{Kh} categorified Jones polynomial. His
construction is a refinement of the Kauffman bracket construction. 
To any classical link diagram $D$ it associates a bigraded chain complex
of abelian groups $C^{i,j}(D)$ with differential 
$d:C^{i,j}(D)\to C^{i+1,j}(D)$ and homology groups $\mathcal H^{i,j}(D)$
such that Reidemeister moves of the diagram induce homotopy equivalences
of the complex. The Khovanov homology $\mathcal H^{i,j}(D)$ categorifies
the Jones polynomial $f_D(A)$ in the sense that
$$f_D(A)=\sum_{i,j}(-1)^{i+j}A^{-2j}\rnk\mathcal H^{i,j}(D). $$

This slightly strange formula is due to the fact that Khovanov used
different normalization of the Jones polynomial, namely
$K(D)(q)=f_D(\frac1{\sqrt{-q}})$. Ranks of the Khovanov homology groups 
are related to $K$ as follows:
$$K(D)(q)=\sum_{i,j}(-1)^iq^j\rnk\mathcal H^{i,j}(D).$$

For any state $s$ of $D$, the construction associates to each connected 
component $C$ of $D_s$ a copy $\cal A_C$ of graded 
free abelian group $\cal A$ with two generators,  $1$ of grading 1 
and $x$ of grading -1. 
To the whole $s$ the construction associates a graded group $C(s)$ which
is the tensor
product of all copies $\cal A_c$ of $\cal A$ associated to the 
connected components of $D_s$ with grading shifted by
$\frac{3w(D)-a(s)+b(s)}2$, where $w(D)$ is the writhe of $D$.

For a state $s$ of diagram $D$ denote $\frac{w(D)-a(s)+b(s)}2$ by
$i(s)$ and consider $C(s)$ as a bigraded group with first grading to be
identically equal to $i(s)$ and the second grading as defined above
(i.e., the grading of tensor product shifted by 
$\frac{3w(D)-a(s)+b(s)}2$). Denote by $\cal C(D)$ the bigraded
group $\sum_sC(s)$. This is the total group of the Khovanov chain
complex.

To define differential, we need to fix an order of all crossings of 
the diagram.
We need to fix also multiplication $m:\cal A\otimes\cal A\to\cal A$
defined by formulas $1\otimes 1\mapsto1$, $1\otimes x\mapsto x$,
$x\otimes 1\mapsto x$ and $x\otimes x\mapsto 0$, and comultiplication
$\Delta:\cal A\to\cal A\otimes\cal A$ defined by formulas $1\mapsto 1\otimes
x+x\otimes1$ and $x\mapsto x\otimes x$.

The differential $C^{i,j}\to C^{i+1,j}$ is defined as the sum of partial
differentials $C^j(s)\to C^j(t)$, where $s$ and $t$ are states with
$i(s)=i$, $i(t)=i+1$ such that $t$ differs from $s$ only by a marker 
at one crossing. Denote this crossing by $x$ and its number by $k$.
Denote by $r$ the number of negative markers of $s$ at crossings with 
numbers greater than $k$.
At $x$ the marker of $s$ is positive, while the marker of $t$ is 
negative. 

Let two connected components of $D_s$ pass near $x$. Denote them by 
$C_1$ and $C_2$. Then there is only one connected component of $D_t$ 
passing near $x$. Denote this component by $C$.
Then $C(s)=\cal A_{C_1}\otimes\cal A_{C_2}\otimes B$ and $C(t)=\cal
A_C\otimes B$ and the partial differential $C(s)\to C(t)$ is defined by
formula $(-1)^rm\otimes id_B$.

Let only one component of $D_s$ pass near $x$. Denote it by $C$.
Then two components of $D_t$ pass near $x$. Denote them by $D_{C_1}$ and 
$D_{C_2}$. Then $C(s)=\cal A_C\otimes B$ and 
$C(t)=\cal A_{C_1}\otimes\cal A_{C_2}\otimes B$ and the partial 
differential $C(s)\to C(t)$ is defined by
formula $(-1)^r\Delta\otimes id_B$.

\subsection{Khovanov Complex of Oriented Signed Chord
Diagram}\label{s5.2}
It is clear that the construction of Khovanov complex sketched above
can be expressed completely in terms of signed chord diagram, as it was
done with the Kauffman bracket in Section \ref{s3.2}. 

If a signed chord diagram is orientable, all properties of its Khovanov
complexes proven in \cite{Kh} or \cite{BN} can be repeated without any
change. In particular, this is a complex. The only property of classical
link diagrams that is used in the proof of this is that change of a
single marker in a state $s$ gives rise to change of the number of 
connected components in $D_s$, and this is a corollary of orientability
of the chord diagram.

Further, the Reidemeister moves of signed chord diagram induce 
homotopy equivalences of the Khovanov complex, provided they 
preserve an orientation. The proof is also borrowed from the classical
case without any change.

According to Theorem \ref{4.7.B}, if alternatable virtual diagrams are
virtualy isotopic then they can be related by an alternatable isotopy.
Therefore, the homotopy type of the Khovanov complex and, in particular,
Khovanov homology groups are invariants of alternatable virtual links.

\subsection{Failure in Non-Orientable Case}\label{s5.3}
If the signed chord diagram is not orientable, there is no problem to
extend the definition of differential, because the partial differentials
along a change of marker destroying orientation of $D_s$ and hence 
preserving the number of components vanish for grading reasons. 

Indeed, 
$C(s)$ is obtained from $\cal A^{|s|}$ by a shift of the second grading by 
$\frac{3w(D)-a(s)+b(s)}2$. The graded rank\footnote{Recall that the graded 
rank of a finitely generated graded group $W=\bigoplus_j W_j$ is the Laurent
polynomial $\sum_j q^j\rnk W_j$.} of $C(s)$ is
$q^{\frac{3w(D)-a(s)+b(s)}2}(q+q^{-1})^{|s|}$. All exponents of this
graded rank are congruent to each other modulo 2. When $a(s)$ and $b(s)$
change by one and $|s|$ does not change,  the graded rank is multiplied
by $q$, and parity of all exponents changes by 1. Thus any homomorphism
preserving the grading is trivial.

However, the homomorphism $\cal C(D)\to\cal C(D)$ obtained in this way
is not a differential: its square is not zero. This can be easily seen
in the very simple example.

The unknot in $\R P^3$ is isotopic to
$\sfig{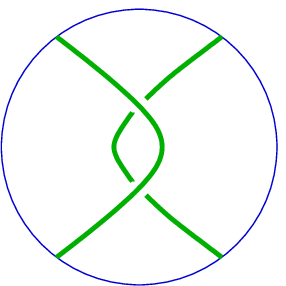}{.7}$.
The corresponding virtual link diagram is shown on the right hand side 
of Figure \ref{A2NA}.

Consider the Khovanov chain groups. There are 4 states 
giving rise to Khovanov chain groups which looks as follows:\\
$$               
\begin{CD}
 & & C\left(\sfig{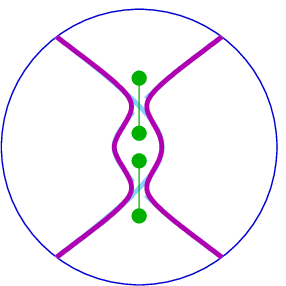}{.6}\right)^{(q+q^{-1})} &&\\
 & \nearrow &                                     & \searrow &  \\
 C\left(\sfig{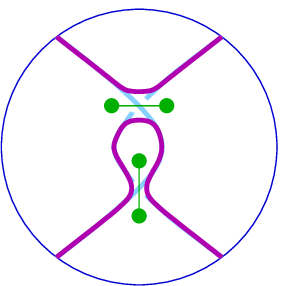}{.6}\right)^{(q^{-2}+1)} &&&&
 C\left(\sfig{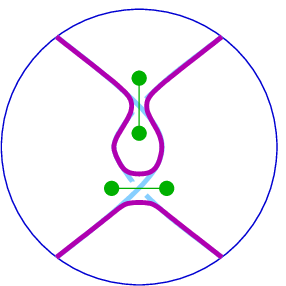}{.6}\right)^{(1+q^2)} \\
& \searrow & & \nearrow &\\
&&   C\left(\sfig{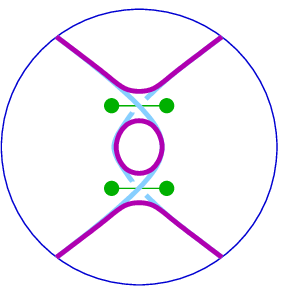}{.6}\right)^{(q^2+2+q^{-2})} && 
\end{CD}
$$
Here the graded ranks are shown as exponents. For the grading reasons
(the partial differentials are of degree 0), the upper arrows should be
zero. Therefore the composition of the bottom two arrows should be zero.
Otherwise the square of the differential would not be zero. 

But it is non-trivial in the component with grading 0. Indeed, the first
of them maps $1$ to $\Delta(1)=1\otimes x+x\otimes 1$, then the second
one maps $1\otimes x$ to $m(1\otimes x)=x$ and $x\otimes 1$ to $x$.
Hence the composition sends $1$ to $2x$.

\end{document}